\documentclass{amsart}

\usepackage{amsmath,amscd,amssymb,amsthm}
\usepackage{enumerate,mathrsfs}
\usepackage{geometry}
\usepackage{hyperref}
\usepackage{url}
\usepackage[all]{xy}

\newcommand{\bibfilename}{/home/micromath/results/b20140922unibib/unibib}
\numberwithin{equation}{section}

\theoremstyle{plain}
\newtheorem{thm_}[equation]{Theorem}
\newtheorem{lemma_}[equation]{Lemma}
\newtheorem{prop_}[equation]{Proposition}
\newtheorem{cor_}[equation]{Corollary}
\newtheorem{eg_}[equation]{Example}
\newtheorem{con_}[equation]{Conjecture}
\newtheorem*{cons_}{Conjecture}

\theoremstyle{definition}
\newtheorem{thmu_}[equation]{Theorem}
\newtheorem*{thmus_}{Theorem}
\newtheorem{propu_}[equation]{Proposition}
\newtheorem*{propus_}{Proposition}
\newtheorem{coru_}[equation]{Corollary}
\newtheorem{lemu_}[equation]{Lemma}
\newtheorem*{lemus_}{Lemma}
\newtheorem{egu_}[equation]{Example}
\newtheorem*{egus_}{Example}
\newtheorem{def_}[equation]{Definition}
\newtheorem*{defs_}{Definition}
\newtheorem{rk_}[equation]{Remark}

\newcommand{\thm}[1]{\begin{thm_}#1\end{thm_}}

\newcommand{\lemm}[1]{\begin{lemma_}#1\end{lemma_}}
\newcommand{\lemu}[1]{\begin{lemu_}#1\end{lemu_}}

\newcommand{\prop}[1]{\begin{prop_}#1\end{prop_}}

\newcommand{\pf}[1]{\begin{proof}#1\end{proof}}

\DeclareMathOperator{\Gal}{Gal}
\DeclareMathOperator{\Hom}{Hom}

\DeclareMathOperator{\Spec}{Spec}

\DeclareMathOperator{\im}{im}

\DeclareMathOperator{\Br}{Br}
\DeclareMathOperator{\Div}{Div}
\DeclareMathOperator{\Pic}{Pic}

\DeclareMathOperator{\coker}{coker}
\DeclareMathOperator{\cone}{cone}
\DeclareMathOperator{\KD}{KD}
\DeclareMathOperator{\U}{U}

\newcommand{\AAb}{\mathbf{Ab}}
\newcommand{\MMod}{\mathbf{Mod}}
\newcommand{\bA}{\mathbb A}%

\newcommand{\HH}{\mathbb H}

\newcommand{\QQ}{\mathbb Q}
\newcommand{\RR}{\mathbb R}
\newcommand{\ZZ}{\mathbb Z}
\newcommand{\bfD}{\mathbf D}%
\newcommand{\bfG}{\mathbf G}%
\newcommand{\bfH}{\mathbf H}%

\newcommand{\bfR}{\mathbf R}
\newcommand{\bfS}{\mathbf S}

\newcommand{\g}{\mathfrak g}

\newcommand{\sF}{\mathscr F}

\newcommand{\cT}{\mathcal T}
\newcommand{\tm}{\times}%
\newcommand{\otm}{\otimes}
\newcommand{\ol}{\overline}
\newcommand{\ra}{\rightarrow}

\newcommand{\Ra}{\Rightarrow}
\newcommand{\lra}{\longrightarrow}

\newcommand{\hra}{\hookrightarrow}

\newcommand{\os}[2]{\overset{#1}{#2}}
\newcommand{\eq}[1]{\begin{equation}#1\end{equation}}
\newcommand{\eqn}[1]{\begin{equation*}#1\end{equation*}}
\newcommand{\ga}[1]{\begin{gather}#1\end{gather}}
\newcommand{\gan}[1]{\begin{gather*}#1\end{gather*}}

\newcommand{\enmt}[1]{\begin{enumerate}#1\end{enumerate}}

\newcommand{\cminus}[1]{#1^c\setminus#1}
\newcommand{\divminus}[1]{\Div_{\cminus{\ol #1}}\ol #1^c}

\begin{document}
\title[Brauer group and the Brauer-Manin set of a product]{A note on the Brauer
 group and the Brauer-Manin set of a product}
\author[C. Lv]{Chang Lv}
\address{State Key Laboratory of Information Security\\
Institute of Information Engineering\\
Chinese Academy of Sciences\\
Beijing 100093, P.R. China}
\email{lvchang@amss.ac.cn}
\date{\today}
\thanks{This work was supported by
 National Natural Science Foundation of China (Grant No. 11701552).
}
\begin{abstract}
We generalize the results of Skorobogatov and Zarhin considering  the commutativity
 of Brauer groups (and Brauer-Manin sets) with taking product of two varieties,
 by  relaxing the condition that varieties are projective.
\end{abstract}
\maketitle

\section{Introduction}\label{sec_intro}
Let $X$  be a variety over a number field $k$, $\Br X=H^2_{\text{\'et}}(X, \bfG_m)$
 the cohomological Brauer-Grothendieck group
 \cite{grothendieck95braueri, grothendieck95brauerii, grothendieck95braueriii},
 $X(\bA_k)^{\Br X}$ the Brauer-Manin set \cite{torsor}.
Consider the product of two varieties $X\tm Y$, and it is natural to ask about
 the commutativity
 of Brauer groups (and Brauer-Manin sets) with taking product.

Skorobogatov and Zarhin \cite{sz14product} investigated
 the cokernel of the natural map $\Br X\oplus\Br Y\ra \Br(X\tm Y)$ and the relation
 between $X(\bA_k)^{\Br X}\tm Y(\bA_k)^{\Br Y}$ and $(X\tm Y)(\bA_k)^{\Br(X\tm Y)}$, namely,
 if $X$ and $Y$ are smooth, projective, geometrically integral,
 then  $\coker(\Br X\oplus\Br Y\ra \Br(X\tm Y))$ is
 finite and   $X(\bA_k)^{\Br X}\tm Y(\bA_k)^{\Br Y}=(X\tm Y)(\bA_k)^{\Br(X\tm Y)}$.

The aim of this note is to relax the projectivity constraint on $X$ and $Y$.
Let $k$  be a field  finitely generated over $\QQ$, with fixed separable closure $\ol k$
 and $\g_k=\Gal(\ol k/k)$.
Let $X$ and $Y$ be smooth quasi-projective
 geometrically integral varieties over $k$ whose base changes to $\ol k$ are denoted by
 $\ol X$ and $\ol Y$.
Then in Proposition \ref{prop_cok_barbrgk} we show that
 $\coker((\Br\ol X)^{\g_k}\oplus (\Br\ol Y)^{\g_k}\ra \Br(\ol X\tm\ol Y)^{\g_k})$
 is finite, which generalizes \cite[Thm. A]{sz14product}.

Assume that $(X\tm Y)(k)\neq\emptyset$ or $H^3(k, \ol k^\tm)=0$.
Then in Theorem \ref{thm_cok_br} we show that
 $\coker(\Br X\oplus\Br Y\ra \Br(X\tm Y))$ is finite,
 generalizing  the first main result Thm. B of \cite{sz14product}.

If $k$  is a number field and  $X$ and $Y$ are smooth geometrically integral varieties
 over $k$, then in Theorem \ref{thm_bmset} we show that
 $(X\tm Y)(\bA_k)^{\Br(X\tm Y)} = X(\bA_k)^{\Br X}\tm Y(\bA_k)^{\Br Y}$,
 generalizing   the second main result Thm. C of \cite{sz14product}.
As mentioned in \cite{sz14product}, by using \'etale homotopy of Artin and Mazur,
 Harpaz and Schlank \cite[Cor. 1.3]{hs13homotopy} proved a statement similar to
 Theorem \ref{thm_bmset} where the Brauer-Manin set is replaced by the \'etale Brauer-Manin
 set and varieties do not need to be proper.

For Theorem \ref{thm_cok_br}, the idea of proof is to compare the various cohomological
 groups of varieties with the ones of their smooth compactifications, which uses a result of
 by Colliot-Th{\'e}l{\`e}ne and Skorobogatov \cite{ctsk13descente}, implying that
 $(\Br\ol X^c)^{\g_k}\ra (\Br\ol X)^{\g_k}$ has finite cokernel,
 where $\ol X^c$ is a smooth compactification of $\ol X$.
Since  the compactifications
 are projective, one can then make use of  the original result of \cite{sz14product}.
The proof of  Theorem \ref{thm_bmset} relies on the original results and the generalizable
 construction \cite{sz14product} of a certain homomorphism using
 universal torsors of $n$-torsion,
 and a modified version (Lemma \ref{lem_H2mun}) of a lemma of Cao \cite{cao18sous}
 considering  the relation between
 $H^2_{\text{\'et}}(X,\mu_n)\oplus H^2_{\text{\'et}}(Y,\mu_n)$
 and $H^2_{\text{\'et}}(X\tm Y,\mu_n)$.

\smallskip
\noindent \textbf{Notations.}
Let $k$ be a field of characteristic zero and
 fix a separable closure $\ol k$. Let $\g_k=\Gal(\ol k/k)$.
A $k$-variety $X$ is a separated $k$-scheme of finite type.
For any field $K$ containing $k$, denote by $X_K$  the base change
 $X\tm_{\Spec k}\Spec K$. We also write  $\ol X$ for $X_{\ol k}$.
For two varieties $X$ and $Y$, we write
 $X\tm Y$ (resp. $\ol X\tm\ol Y$) for $X\tm_{\Spec k} Y$
 (resp. $\ol X\tm_{\Spec\ol k}\ol Y$).
Denote by $p_X$ and $p_Y$ the two projections from
 $X\tm Y$ to $X$ and $Y$.

Since $k$ is of characteristic zero, if $X$ is smooth, by Hironaka's theorem, the smooth
 compactification  of $\ol X$ exists.
Let  $\ol X^c$  be a  smooth compactification  of $\ol X$. Then we have
 a dense open immersion $\ol X\hra \ol X^c$.
If  $X$ is also quasi-projective, $\ol X^c$  can be made projective.

Denote by $\bfS(X)$, $\AAb$ and $\g_k$-$\MMod$ for the categories of \'etale sheaves on $X$,
 abelian groups and discrete $\g_k$-modules.
Let $\bfD^+(X)$ and $\bfD^+(\AAb)$  be the corresponding derived categories
 of complexes bounded below.
Since $\bfS(X)$ has enough injectives, the right derived functor
 $\bfR F:\bfD^+(X)\ra\bfD^+(\AAb)$ of any additive functor $F:\bfS(X)\ra\AAb$ exists.
We write $\RR^i F=H^i\bfR F$  for the $i$-th hypercohomology functor.
Let $\bfH(X,-)=\bfR\Gamma(X_{\text{\'et}},-)$ and
 $\HH^i(X,-)=\RR^i\Gamma(X_{\text{\'et}},-)$, where $X_{\text{\'et}}$ is the small \'etale
 site over $X$.
Let $\bfS(k)=\bfS(\Spec k)\cong \g_k$-$\MMod$, $\bfD^+(k)=\bfD^+(\Spec k)$.
For $\sF\in\bfS(X)$, let $H^i(X,\sF)=\HH^i(X,\sF)$,
  and for $\sF\in\bfS(k)$,  $H^i(k, \sF)=\HH^i(\Spec k, \sF)$.

If $n$ is a positive number and $A$ an abelian group, write $A_n$ and $A/n$ for
 the kernel and cokernel of the homomorphism $A\os{\tm n}{\ra} A$.

\section{The finiteness of the cokernel of the Brauer groups} \label{sec_cok}
For a contravariant functor $F$ from  the category of $k$-varieties to $\AAb$,
let $X$ and $Y$ be two $k$-varieties with compactifications $\ol X^c$ and $\ol Y^c$.
Define
\gan{
\coker(F)=\coker( F(X)\oplus F(Y)\ra F(X\tm Y) ), \\
\ol\coker(F)=\coker( F(\ol X)\oplus F(\ol Y)\ra F(\ol X\tm \ol Y) ), \\
\ol\coker^c(F)=\coker( F(\ol X^c)\oplus F(\ol Y^c)\ra F(\ol X^c\tm \ol Y^c) )
}
 where the maps are induced by  $(p_X^*, p_Y^*)$. We use similar notation for $\ker$.

Define $\Br_1 X$ and $\Br_2 X$ to be the kernel and cokernel of the natural map
 $\Br X\ra (\Br\ol X)^{\g_k}$.

We first consider $\ol\coker(\Br^{\g_k})$.
\lemm{\label{lem_Hn_inj}
Let   $X$  and $Y$ be two varieties over $k$ and  $G$  a \'etale sheaf defined by a
 commutative $k$-group scheme.
Then for $n\ge1$,  we have $\ol\ker(H^n(-, \ol G))=0$, that is,
 the natural map  induced by  $(p_X^*, p_Y^*)$
\eqn{
H^n(\ol X, \ol G)\oplus H^n(\ol Y, \ol G)\ra H^n(\ol X\tm\ol Y, \ol G)
}
 is an injective homomorphism of $\g^k$-modules.
}
\pf{
See \cite[Prop. 1.5 (iv)]{sz14product} where varieties are assumed to be projective.
One checks that it holds for general varieties.
}

\prop{\label{prop_cok_barbrgk}
Let $k$  be a field finitely generated over $\QQ$. Let $X$ and $Y$ be smooth quasi-projective
 geometrically integral varieties over $k$.
Then the group $\ol \coker(\Br^{\g_k})$, which is the cokernel of
\eqn{
(\Br\ol X)^{\g_k}\oplus (\Br\ol Y)^{\g_k}\ra \Br(\ol X\tm\ol Y)^{\g_k},
}
 is finite.
}
\pf{
Since in the projective case, $\ol \coker^c(\Br^{\g_k})$  is already finite
 \cite[Thm. A]{sz14product}, the result follows from  Lemmas
 \ref{lem_cok_br_fin} \eqref{it_cok_brcbrgk_brgk_fin} and \ref{lem_alpha_beta} below,
 which imply that  $\ol \coker^c(\Br^{\g_k})$  is finite if and only if
 $\ol \coker(\Br^{\g_k})$   is.
}

\lemm{\label{lem_cok_br_fin}
Let $k$  be a field finitely generated over $\QQ$ and
 $X$ a smooth quasi-projective geometrically integral variety over $k$.
Then we have
\enmt{[\upshape (i)]
\item \label{it_cok_brcbrgk_brgk_fin}
 the  natural  homomorphism  $\alpha(X): (\Br\ol X^c)^{\g_k}\ra (\Br\ol X)^{\g_k}$
 induced by $\ol X\hra \ol X^c$ is injective and has finite cokernel,
\item \label{it_br2_fin}
 the  natural  homomorphism  $\Br X\ra (\Br\ol X)^{\g_k}$ has finite cokernel,
 i.e., $\Br_2X$ is finite.
}}
\pf{
The second statement of the lemma is \cite[Thm. 6.2 (i)]{ctsk13descente},
 whose proof, based on the proper base change theorem and Weil conjecture, contains
 the  first statement of the lemma (p. 164, l. 4-5, \emph{loc. cit.}).
}

\lemm{\label{lem_alpha_beta}
With notations as Proposition \ref{prop_cok_barbrgk},
 let $\beta: \ol\coker^c(\Br^{\g_k})\ra \ol\coker(\Br^{\g_k})$ be the natural homomorphism
 induced by  $\ol X\hra \ol X^c$ and $\ol Y\hra \ol Y^c$.
Then we have the exact sequence
\eqn{
0\ra \ker\beta\ra \coker\alpha(X)\oplus\coker\alpha(Y)\ra \coker\alpha(X\tm Y)\ra
 \coker\beta\ra 0,
}
 where $\alpha$ is defined in Lemma \ref{lem_cok_br_fin}
 \eqref{it_cok_brcbrgk_brgk_fin}.
}
\pf{
By Lemma \ref{lem_cok_br_fin} \eqref{it_cok_brcbrgk_brgk_fin},
 there exists a functorial short exact sequence
\eq{\label{eq_barbrgk_compare}
0\ra (\Br\ol X^c)^{\g_k}\os{\alpha(X)}{\ra} (\Br\ol X)^{\g_k}\ra \coker\alpha(X)\ra 0.
}
Next, by Lemma \ref{lem_Hn_inj} we have the functorial short exact sequence
\eq{\label{eq_cok_barbrgk}
0\ra (\Br\ol X)^{\g_k}\oplus (\Br\ol Y)^{\g_k}\ra  (\Br\ol X\tm\ol Y)^{\g_k}\ra
 \ol\coker(\Br^{\g_k})\ra 0,
}
 which fits into  the commutative diagram with exact rows
\eqn{\xymatrix{
0\ar[r] &(\Br\ol X^c)^{\g_k}\oplus (\Br\ol Y^c)^{\g_k}\ar[r]
 \ar[d]^-{\alpha(X)\oplus\alpha(Y)}
 &(\Br\ol X^c\tm\ol Y^c)^{\g_k}\ar[r]\ar[d]^-{\alpha(X\tm Y)}
 &\ol\coker^c(\Br^{\g_k})\ar[r]\ar[d]^-{\beta} &0 \\
0\ar[r] &(\Br\ol X)^{\g_k}\oplus (\Br\ol Y)^{\g_k}\ar[r]
 &(\Br\ol X\tm\ol Y)^{\g_k}\ar[r] &\ol\coker(\Br^{\g_k})\ar[r] &0
}}
The desired exact sequence then follows  from
 \eqref{eq_barbrgk_compare} and the snake lemma.
}

Next, we investigate $\coker(\Br_1)$.
\prop{\label{prop_cok_br1}
Let $k$  be a field finitely generated over $\QQ$.
Let $X$ and $Y$ be smooth quasi-projective
 geometrically integral varieties over $k$.
Assume that $(X\tm Y)(k)\neq\emptyset$ or $H^3(k, \ol k^\tm)=0$.
Then $\coker(\Br_1)$, the cokernel of the natural map
\eqn{
\Br_1 X\oplus\Br_1 Y\ra \Br_1(X\tm Y),
}
 is finite.
}
For the proof, it is enough to use Lemmas \ref{lem_cokbr1_iso_cokH1kKD'},
 \ref{lem_tri_KD'_cokPic} and \ref{lem_cokcPic_iso_cokPic} below.

For an arbitrary $k$-variety $X$, let $p: X\ra \Spec k$ be the structure morphism.
Following Harari and Skorobogatov \cite{hs13descent},
 define $\KD(X)=(\tau_{\le1}\bfR p_*\bfG_{m,X})[1]$ and
\eqn{
\KD'(X)=\cone( \bfG_{m,k}\ra \tau_{\le1}\bfR p_*\bfG_{m,X} )[1].
}
Thus we have the distinguished triangle in $\bfD^+(k)$
\eq{\label{eq_tri_KD'}
\bfG_m[1]\ra \KD(X)\ra \KD'(X)\ra.
}

\lemm{\label{lem_cokbr1_iso_cokH1kKD'}
Let $k$  be a field of characteristic zero. Let $X$ and $Y$ be  two $k$-varieties.
Assume that $(X\tm Y)(k)\neq\emptyset$ or $H^3(k, \ol k^\tm)=0$.
Then we have a natural isomorphism $\coker(\Br_1)\os{\sim}{\ra} \coker(\HH^1(k, \KD'))$.
}
\pf{
By applying $\bfH(k,-)$ to the distinguished triangle \eqref{eq_tri_KD'} and
 taking cohomology,
the exact sequence  \cite[(8.5)]{hs13descent} extends to 
\eqn{
\Br k\ra\Br_1 X\ra\HH^1(k, \KD'(X))\ra H^3(k, \ol k^\tm)\os{e}{\ra}
 \HH^2(k, \KD(X)). 
}
We have an obvious commutative diagram
\eqn{\xymatrix{
\bfG_m[1]\ar[r]\ar[rd] &\KD(X)\ar[d] \\
 &\bfR p_*\bfG_m[1]
}}
On applying $\HH^2(k,-)$  and using $\HH(k,\bfR p_*-)=\HH(X,-)$ we obtain
\eqn{\xymatrix{
H^3(k,\ol k^\tm)\ar[r]^-{e}\ar[rd]^-{p^*} &\HH^2(k, \KD(X))\ar[d] \\
 &H^3(X, \bfG_m)
}}
 where  $p^*$ is an injection by the assumption that
 $(X\tm Y)(k)\neq\emptyset$ or $H^3(k, \ol k^\tm)=0$
 (see the comments after (21) in \cite{sz14product}),
 and so is the horizontal arrow $e$.
Thus by functoriality we have the commutative diagram with exact rows
\eqn{\xymatrix@-2mm{
\Br k\oplus \Br k \ar[r]\ar@{->>}[d] &\Br_1 X\oplus \Br_1 Y\ar[r]\ar[d]
 &\HH^1(k, \KD'(X))\oplus \HH^1(k, \KD'(Y))\ar[r]\ar[d] &0\\
\Br k\ar[r] &\Br_1(X\tm Y)\ar[r] &\HH^1(k, \KD'(X\tm Y)\ar[r] &0
}}
Clearly the left vertical arrow is surjective.
Thus the proof is complete by using the snake lemma.
}

\lemm{\label{lem_tri_KD'_cokPic}
Let $k$  be a field of characteristic zero. Let $X$ and $Y$ be  geometrically
 integral $k$-varieties.
Then we have a  distinguished triangle in $\bfD^+(k)$
\eqn{
\KD'(X)\oplus \KD'(Y)\ra \KD'(X\tm Y)\ra \ol\coker(\Pic)\ra.
}
In particular,  we have a natural injection
 $\coker(\HH^1(k, \KD'))\hra H^1(k, \ol\coker(\Pic))$.
}
\pf{
Let $\U(X)=\cone(\bfG_m\ra p_*\bfG_m)$.
We have  the obvious distinguished triangles
\gan{
\bfG_m[1]\ra p_*\bfG_m[1]\ra  \U(X)[1]\ra, \\
p_*\bfG_m[1]\ra \KD(X)\ra (\tau_{[1]}\bfR p_*\bfG_m)[1]\ra, \\
\bfG_m[1]\ra \KD(X)\ra \KD'(X)\ra.
}
Note $(\tau_{[1]}\bfR p_*\bfG_m)[1]=R^1 p_*\bfG_m=\Pic\ol X$.
Then by the octahedron axiom of triangulated categories,
 we have the distinguished triangle
\eqn{
\U(X)[1]\ra \KD'(X)\ra \Pic\ol X\ra.
}
Thus  by functoriality we have the commutative diagram of distinguished triangles
\eqn{\xymatrix{
\U(X)[1]\oplus\U(Y)[1]\ar[r]\ar[d] &\KD'(X)\oplus\KD'(Y)\ar[r]\ar[d]^-{g}
 &\Pic\ol X\oplus\Pic\ol Y\ar[r]\ar[d] & \\
\U(X\tm Y)[1]\ar[r] &\KD'(X\tm Y)\ar[r]\ar[d] &\Pic(\ol X\tm \ol Y)\ar[r]\ar[d] & \\
 &\cone(g)\ar[r]\ar[d] &\ol\coker(\Pic)\ar[d] \\
 & &
}}
Next, by Rosenlicht lemma (c.f. \cite[Lem. 6.5 (ii)]{sansuc81groupe}),
 we obtain a natural isomorphism
\eq{\label{eq_U}
 \U(X\tm Y)\os{\sim}{\ra}\U(X)\oplus \U(Y).
}
Thus  the desired distinguished triangle is obtain using the octahedron axiom.

Applying $\bfH(k,-)$ we have the long exact sequence
\gan{
\dots\ra \HH^i(k, \KD'(X))\oplus\HH^i(k, \KD'(X))\ra \HH^i(k, \KD'(X\tm Y))\ra \\
 H^i(k, \ol\coker(\Pic))\ra \HH^{i+1}(k, \KD'(X))\oplus\HH^{i+1}(k, \KD'(X))\ra \dots.
}
Taking $i=1$ we obtain the desired injection. The proof is complete.
}

\lemm{\label{lem_cokcPic_iso_cokPic}
Let $k$  be a field finitely generated over $\QQ$, and let
 $X$ and $Y$ be smooth quasi-projective geometrically integral varieties over $k$.
Then we have a natural isomorphism of $\g_k$-modules
 $\ol\coker^c(\Pic)\os{\sim}{\ra}\ol\coker(\Pic)$
 induced by  $\ol X\hra \ol X^c$ and $\ol Y\hra \ol Y^c$.
In particular, as an abelian group,  $\ol\coker(\Pic)$ is  finitely generated and
 torsion-free, and  $H^1(k, \ol\coker(\Pic))$ is finite.
}
\pf{
Let $\divminus{X}=\bigoplus_{x\in (\ol X^c)^{(1)}\cap(\cminus{\ol X})}\ZZ$ be the group of
 divisors of $\ol X^c$ outside $\ol X$,
 where $(\ol X^c)^{(1)}$ is the set of codimesion $1$-point in
 $\ol X^c$.
Since $X$ is smooth,  we  have  a natural exact sequence
 (c.f. \cite[(4.18), (4.19), pp. 72]{torsor})
\eqn{
0\ra \U(X)\ra \divminus{X}\ra \Pic\ol X^c\ra \Pic\ol X\ra 0
}
 which, by functoriality and Lemma \ref{lem_Hn_inj},
 fits into the following exact commutative diagram of $\g_k$-modules
\eqn{\xymatrix@-4mm{
 & & &0\ar[d] &0\ar[d] \\
0\ar[r] &\U(X)\oplus\U(Y)\ar[r]\ar[d]^{\wr} &\divminus{X}\oplus\divminus{Y}\ar[r]
 \ar[d]^-{\wr} &\Pic\ol X^c\oplus \Pic\ol Y^c\ar[r]\ar[d]
 &\Pic\ol X\oplus\Pic\ol Y\ar[r]\ar[d] &0 \\
0\ar[r] &\U(X\tm Y)\ar[r] &\Div_{\ol X^c\tm\ol Y^c\setminus\ol X\tm\ol Y}\ol X^c\tm\ol Y^c
 \ar[r] &\Pic(\ol X^c\tm\ol Y^c)\ar[r]\ar[d] &\Pic(\ol X\tm\ol Y)\ar[r]\ar[d] &0 \\
 & & &\ol\coker^c(\Pic)\ar[r]\ar[d] &\ol\coker(\Pic)\ar[d] \\
 & & &0 &0
}}
 where the first vertical arrow is an isomorphism by \eqref{eq_U}
 and so is the second one since by definition,
 $(\ol X^c\tm\ol Y^c)^{(1)}\cap (\ol X^c\tm\ol Y^c\setminus \ol X\tm\ol Y) =
 ((\ol X^c)^{(1)}\cap(\cminus{\ol X}))\sqcup ((\ol Y^c)^{(1)}\cap(\cminus{\ol Y}))$.
Then the desired isomorphism  follows from the diagram and the snake lemma.

Next, denote by $A$ the Picard variety of $X^c$ and $B^t$ the Albanese variety of $Y^c$.
Then the abelian group $\ol\coker^c(\Pic)=\Hom(\ol{B^t},\ol A)$ is finitely generated
 and torsion-free \cite[Prop. 1.7 and Rem. 1.9]{sz14product}, and so is $\ol\coker(\Pic)$.
It follows that  $H^1(k, \ol\coker(\Pic))$ is finite. The proof is complete.
}

Now we obtain the desired finiteness for $\coker(\Br)$.
\thm{\label{thm_cok_br}
Let $k$  be a field finitely generated over $\QQ$. Let $X$ and $Y$ be smooth quasi-projective
 geometrically integral varieties over $k$.
Assume that $(X\tm Y)(k)\neq\emptyset$ or $H^3(k, \ol k^\tm)=0$.
Then the group $\coker(\Br)$, which is the cokernel of the natural map
\eqn{
\Br X\oplus\Br Y\ra \Br(X\tm Y),
}
 is finite.
}
\pf{
Consider  the following commutative diagram with exact rows
\eqn{\xymatrix@-4mm{
0\ar[r] &\Br_1X\oplus\Br_1Y\ar[r]\ar[d] &\Br X\oplus\Br Y\ar[r]
 \ar[d] &(\Br\ol X)^{\g_k}\oplus(\Br\ol Y)^{\g_k}
 \ar[r]\ar[d] &\Br_2X\oplus\Br_2Y\ar[r]\ar[d] &0\\
0\ar[r] &\Br_1(X\tm Y)\ar[r] &\Br(X\tm Y)\ar[r]
 &\Br(\ol X\tm\ol Y)^{\g_k}\ar[r] &\Br_2(X\tm Y)\ar[r] &0
}}
The idea is to break rows into short exact sequences and use the snake lemma.
By Lemma  \ref{lem_cok_br_fin} \eqref{it_br2_fin}, $\ker(\Br_2)$ is finite.
Proposition   \ref{prop_cok_barbrgk} implies $\ol\coker(\Br^{\g_k})$ is finite.
Also by Proposition \ref{prop_cok_br1}, $\coker(\Br_1)$ is finite.
Finally  we conclude that $\coker(\Br)$ is finite. The proof is complete.
}

\section{The Brauer-Manin set  of a product} \label{sec_bmset}
\thm{\label{thm_bmset}
Let $k$  be a number field. Let $X$ and $Y$ be smooth geometrically
 integral varieties over $k$. Then
\eqn{
(X\tm Y)(\bA_k)^{\Br(X\tm Y)} = X(\bA_k)^{\Br X}\tm Y(\bA_k)^{\Br Y}.
}}
\pf{
Note that the Brauer group of a  smooth variety is torsion
 \cite[Prop. 1.4]{grothendieck95brauerii}, so we only need to show that
 \cite[Lems. 5.2, 5.3]{sz14product} are still available for non-projective varieties.
Actually, the proof of Lemma 5.3 did not make any use of projective property.
It suffices to prove Lemma \ref{lem_H2mun} below, a variant of Lemma 5.2,
 \emph{loc. cit.}, which does \emph{not} need the projective assumption on varieties.
}
Let $S_X$ be the $k$-group of multiplicative type dual to $H^1(\ol X,\mu_n)$, that is,
\eqn{
S_X^*=\Hom_{\text{$\ol k$-gps}}(S_X,\bfG_m)=H^1(\ol X,\mu_n)
}
 and the same for $S_Y$.
Let  $\cT_X$ be a \emph{universal torsor of $n$-torsion} for $X$, which is a torsor
 over $X$ under $S_X$. It can be defined for a general smooth geometrically integral
 variety over a field of characteristic zero \cite[Def. 2.1]{cao18sous}, and
 exists if $X(\bA_k)^{\Br_1 X}\neq\emptyset$ \cite[Cor. 8.17]{hs13descent}.
Then there is a  homomorphism
\eqn{
\epsilon: \Hom_k(S_X,S_Y^*)\ra H^2(X\tm Y,\mu_n): \phi\mapsto\phi_*[\cT_X]\cup[\cT_Y],
}
 which is originally defined by Skorobogatov and Zarhin \cite{sz14product} and
 can be generalized for  arbitrary smooth geometrically integral varieties.
By \cite[Prop. 2.2]{cao18sous}, we have $[\ol{\cT_X}]=[\cT_{\ol X}]$ in
 $H^1(\ol X, S_X)$ and there is a natural isomorphism
\eqn{
\tau_X: S_X^*=\Hom_{\ol k}(S_X, \mu_n) \os{\sim}{\ra} H^1(\ol X, \mu_n):
 \phi\mapsto\phi_*([\cT_{\ol X}]).
}
Twisting by $\ZZ/n(-1)=\mathscr{H}om(\mu_n, \ZZ/n)$, we obtain the isomorphism
\eqn{
\tau_X(-1): \Hom_{\ol k}(S_X, \ZZ/n) \os{\sim}{\ra} H^1(\ol X, \ZZ/n):
 \psi\mapsto\psi_*([\cT_{\ol X}]).
}
We also obtain that pairing with $[\cT_{\ol X}]$ induces the identity on $H^1(\ol X, \mu_n)$
 and $H^1(\ol X, \ZZ/n)$.
Similar properties hold if we replace $X$  by  $Y$.
Thus  pairing  with $[\cT_{\ol X}]\cup[\cT_{\ol Y}]$ induces the cup product
\eqn{
H^1(\ol X, \ZZ/n)\otm_{\ZZ/n} H^1(\ol Y, \mu_n)\os{\cup}{\lra}H^2(\ol X\tm \ol Y, \mu_n)
}
 and we have the commutative diagram
\eq{\label{eq_epsilon_cup}
\xymatrix{
(\Hom_{\ol k}(S_X,\ZZ/n)\otm_{\ZZ/n} S_Y^*)^{\g_k}\ar@{=}[r]
 \ar[dr]^-{\sim}_-{\tau_X(-1)\otm\tau_Y}
 &\Hom_k(S_X,S_Y^*)\ar[r]^-{\epsilon} &H^2(X\tm Y,\mu_n)\ar[d] \\
 &(H^1(\ol X, \ZZ/n)\otm_{\ZZ/n} H^1(\ol Y, \mu_n))^{\g_k}\ar[r]^-{\cup}
 &H^2(\ol X\tm \ol Y, \mu_n)^{\g_k}
}}
The above discussion can be found in Section 2, \emph{loc. cit.}.

Let us first recall
\lemu{[{\cite[Lemma  5.3]{sz14product}}]
With notations and assumptions as Theorem \ref{thm_bmset},
 for any positive integer $n$, we have
\eqn{
X(\bA_k)^{(\Br_1 X)_n}\tm Y(\bA_k)^{(\Br_1 Y)_n}\subseteq (X\tm Y)(\bA_k)^{\im{\epsilon}}.
}}

\lemm{\label{lem_kunneth}
Assume that $k$ is a field of characteristic zero and $X$, $Y$ are smooth
 geometrically integral $k$-varieties.
Then we have the following \emph{K\"unneth} decompositions of $\g_k$-modules
\ga{
(p_X^*,p_Y^*): H^1(\ol X, \mu_n)\oplus H^1(\ol Y, \mu_n)\os{\sim}{\ra}
 H^1(\ol X\tm\ol Y, \mu_n), \label{eq_H1_dec} \\
(p_X^*, p_X^*\cup p_Y^*, p_Y^*): H^2(\ol X, \mu_n)\oplus H^2(\ol Y, \mu_n)
 \oplus ( H^1(\ol X, \ZZ/n)\otm_{\ZZ/n} H^1(\ol Y, \mu_n) )
 \os{\sim}{\ra} H^2(\ol X\tm\ol Y, \mu_n), \label{eq_H2_dec}
}
 where $\cup$ is the cup product.
}
\pf{
The original proof assuming projectivity is in \cite{sz14product}.
A proof without properness assumption can be found in \cite[Prop. 2.6]{cao18sous}.
}

\lemm{\label{lem_H2mun}
With notations as Theorem \ref{thm_bmset}, suppose that
 $X(\bA_k)^{\Br_1 X}$ and $Y(\bA_k)^{\Br_1 Y}$ are both non-empty.
Then we have a natural surjection
\eq{\label{eq_H2_surj}
\xymatrix@+4mm{
H^2(X,\mu_n)\oplus H^2(Y,\mu_n)\oplus\Hom_k(S_X, S_Y^*)
 \ar@{->>}[r]^-{(p_X^*,p_Y^*,\epsilon)} &H^2(X\tm Y,\mu_n).
}}}
\pf{
We have seen the assumption that
 $X(\bA_k)^{\Br_1 X}$ and $Y(\bA_k)^{\Br_1 Y}$ are both non-empty ensures the
 existence of $\cT_X$ and $\cT_Y$.

It should be noted  that if we assume that $X(k)\neq\emptyset$, fix $u\in X(k)$,
 and define
\eqn{
H^i_u(X,\mu_n)=\ker( H^i(X,\mu_n)\os{u^*}{\ra} H^i(k,\mu_n) ).
}
Then \cite[Cor. 2.7]{cao18sous} asserts that there is an isomorphism
\eq{\label{eq_H2iso}
\xymatrix@+4mm{
H_u^2(X,\mu_n)\oplus H^2(Y,\mu_n)\oplus\Hom_k(S_X, S_Y^*)
 \ar@{->}[r]^-{(p_X^*,p_Y^*,\epsilon)}_-{\sim} &H^2(X\tm Y,\mu_n),
}}
 where $X$, $Y$ are $U$, $V$ there. Thus clearly \eqref{eq_H2_surj} is surjective.

Now it is not necessarily that $X(k)\neq\emptyset$, but by modifying  the proof of
 Cor. 2.7, \emph{loc. cit.}, we can also show  that \eqref{eq_H2_surj} is surjective.
Consider the   Hochschild-Serre spectral sequences
\eq{\label{eq_E2pqX}
E_2^{p,q}(X)=H^p(k, H^q(\ol X, \mu_n))\Ra E^{p+q}(X)=H^{p+q}(X, \mu_n)
}
 for $X$ and $E_2^{p,q}(Y)$ (resp. $E_2^{p,q}(X\tm Y)$) for $Y$ (resp. $X\tm Y$).
Let
\eqn{
\phi_2^{p,q}: E_2^{p,q}(X)\oplus E_2^{p,q}(Y)\ra  E_2^{p,q}(X\tm Y)
}
 be the   morphism of  spectral sequences induced by $(p_X^*, p_Y^*)$.
Let $E_1^2(X)=\ker(E^2(X)\ra E_2^{0,2}(X))$ and  we define $E_1^2(Y)$ and $E_1^2(X\tm Y)$
 in the same manner.

Let $p: X\ra\Spec k$ be the structure morphism.
Applying $\bfH(k,-)$ to the distinguished triangle  
\eqn{
\tau_{[0]}\bfR p_*\mu_n\ra \tau_{\le 1}\bfR p_*\mu_n\ra
 \tau_{[1]}\bfR p_*\mu_n\ra
}
 and taking cohomology, we obtain the exact sequence
\eq{\label{eq_low_ext}
E_2^{2,0}(X)\ra \HH^2(k, \tau_{\le1}\bfR p_*\mu_n)\ra E_2^{1,1}(X)\ra
 E_2^{3,0}(X)\os{f}{\ra} \HH^3(k, \tau_{\le1}\bfR p_*\mu_n)
}
 which extends the low term exact sequence associated to the spectral
 sequence $E_2^{p,q}(X)$ since by a similar argument as in
 \cite[pp. 413-414]{skorobogatov99beyond} we have
 $\HH^2(k, \tau_{\le1}\bfR p_*\mu_n)=E_1^2(X)$.

We have an obvious commutative diagram
\eqn{\xymatrix{
\tau_{[0]}\bfR p_*\mu_n\ar[r]\ar[rd] &\tau_{\le 1}\bfR p_*\mu_n\ar[d] \\
 &\bfR p_*\mu_n
}}
On applying $\HH^3(k,-)$, we obtain
\eqn{\xymatrix{
E_2^{3,0}(X)\ar[r]^-{f}\ar[rd]^-{p^*} &\HH^3(k, \tau_{\le 1}\bfR p_*\mu_n)\ar[d] \\
 &E^3(X)
}}
By assumption $X(\bA_k)\neq\emptyset$, we have the injection
 $p^*: H^3(k,\mu_n)\hra H^3(X,\mu_n)$ by
  \cite[p. 765, l. 6-14]{sz14product}.
Note that we have
\eq{\label{eq_H0_mun}
H^0(\ol X, \mu_n)=\mu_n(\ol k).
}
Hence  $E_2^{3,0}(X)=H^3(k, \mu_n)\ra E^3(X)=H^3(X, \mu_n)$ is injective.
It follows that $f$ is also an injection and  \eqref{eq_low_ext} becomes
\eqn{
E_2^{2,0}(X)\ra E_1^2(X)\ra E_2^{1,1}(X)\ra 0
}
 which by functoriality fits into the commutative diagram with exact rows and vertical
 arrows induced by $\phi_2^{p,q}$
\eqn{\xymatrix@-2mm{
E_2^{2,0}(X)\oplus E_2^{2,0}(Y)\ar[r]\ar@{->>}[d]^-{\phi_2^{2,0}}
 &E_1^2(X)\oplus E_1^2(Y)\ar[r]\ar[d]
 &E_2^{1,1}(X)\oplus E_2^{1,1}(Y)\ar[r]\ar[d]_-{\wr}^-{\phi_2^{1,1}} &0 \\
E_2^{2,0}(X\tm Y)\ar[r] &E_1^2(X\tm Y)\ar[r] &E_2^{1,1}(X\tm Y)\ar[r] &0
}}
Since for any $p$, $\phi_2^{p, 0}$ is clearly surjective by \eqref{eq_H0_mun}
 and $\phi_2^{p,1}$ is
 an isomorphism by \eqref{eq_H1_dec},  the left vertical arrow is surjective
 and the right one is an isomorphism.
It follows by the five lemma that the middle one is surjective.

Using the injectivity of $E_2^{3,0}(X)\ra E^3(X)$ again,
 a standard computation of the spectral sequence \eqref{eq_E2pqX}
 yields the exact sequence
\eqn{
0\ra E_1^2(X)\ra E^2(X)\ra E_2^{0,2}(X)\ra E_2^{2,1}(X).
}
Thus we have the commutative diagram with exact rows
\eqn{\xymatrix@-2mm{
0\ar[r] &E_1^2(X)\oplus E_1^2(Y)\ar[r]\ar@{->>}[d] &E^2(X)\oplus E^2(Y)\ar[r]\ar[d]
 &E_2^{0,2}(X)\oplus E_2^{0,2}(Y)\ar[r]\ar[d]^-{\phi_2^{0,2}}
 &E_2^{2,1}(X)\oplus E_2^{2,1}(Y)\ar[d]_-{\wr}^-{\phi_2^{2,1}} \\
0\ar[r] &E_1^2(X\tm Y)\ar[r] &E^2(X\tm Y)\ar[r] &E_2^{0,2}(X\tm Y)\ar[r]
 &E_2^{2,1}(X\tm Y)
}}
 where the vertical arrows are induced by $\phi_2^{p,q}$, of which the first one is
 a surjection and the last one isomorphism, by the previous discussion.
A diagram chasing yields the exact sequence
\eq{\label{eq_E2_phi}
E^2(X)\oplus E^2(Y)\os{(p_X^*,p_Y^*)}{\lra} E^2(X\tm Y)\ra \coker\phi_2^{0,2}.
}
Then by \eqref{eq_H2_dec} and commutative diagram \eqref{eq_epsilon_cup},
 one shows that the composition
\eqn{
\Hom_k(S_X,S_Y^*)\os{\epsilon}{\ra}  H^2(X\tm Y,\mu_n)\ra \coker \phi_2^{0,2}
}
 is an isomorphism. Along with \eqref{eq_E2_phi} we obtain the desired result that
 \eqref{eq_H2_surj} is surjective.
}

\section*{Acknowledgment}
The author would like to thank  the referees for valuable suggestions,
 and Dasheng Wei, Junchao Shentu and Jiangxue Fang for  helpful discussions.

\bibliography{\bibfilename}
\bibliographystyle{amsplain}
\end{document}